%% file: paper.tex
\documentclass[letterpaper,10pt,conference]{ieeeconf}

\IEEEoverridecommandlockouts  
\usepackage{graphicx}
\usepackage{dsfont} 
\usepackage{amsmath}
\usepackage{amssymb}
\usepackage{url}
\usepackage{algpseudocode}
\usepackage[utf8]{inputenc} 
\usepackage[format=hang,singlelinecheck=false,caption=false,font=footnotesize]{subfig}
\usepackage{changes}

\graphicspath{{./figures/}}

\let\originalleft\left
\let\originalright\right
\DeclareRobustCommand{\left}{\mathopen{}\mathclose\bgroup\originalleft}
\DeclareRobustCommand{\right}{\aftergroup\egroup\originalright}

\newcommand{\R}{\mathbb{R}}
\newcommand{\N}{\mathbb{N}}
\renewcommand{\epsilon}{\varepsilon}
\renewcommand{\phi}{\varphi}
\newcommand{\prno}{\mathit{Pr}}
\newcommand{\reno}{\mathit{Re}}
\newcommand{\derivoper}{\mathrm{d}}

\DeclareMathOperator{\sint}{int}

\newcommand{\norm}[1]{\left\lVert #1\right\rVert}

\newcommand{\set}[1]{\left\{ #1\right\}}
\newcommand{\setb}[1]{\bigl\{ #1\bigr\}}

\let\p=\paren

\let\pb=\parenb

\newcommand{\sqparen}[1]{\left[ #1\right]}
\newcommand{\sqparenb}[1]{\bigl[ #1\bigr]}

\newcommand{\dprod}[2]{\left\langle #1, #2\right\rangle}
\newcommand{\dprodb}[2]{\bigl\langle #1, #2\bigr\rangle}

\let\smats=\sqmatrixs

\newcommand{\pderiv}[3][]{\frac{\partial^{#1} #2}{\partial #3^{#1}}}
\newcommand{\diff}[1]{\derivoper #1}

\newcommand{\diver}[2][]{\nabla_{\!#1} \cdot #2}

\newcommand{\grad}[2][]{\nabla_{\!#1} #2}
\newcommand{\lapl}[2][]{\triangle_{#1} #2}
\newcommand{\bdry}[1]{\partial #1}
\newcommand{\unit}[1]{\sqparen{\rm #1}}

\newcommand{\uf}[1]{\textbf{DEPRECATED}}
\newcommand{\bfn}[1]{\widehat{#1}}

\title{\LARGE \bf%
  Zoned HVAC Control via PDE-Constrained Optimization%
}

\author{%
  Runxin He \and
  Humberto Gonzalez%
  \thanks{%
    The authors are with the Department of Electrical \& Systems Engineering, Washington University in St.\ Louis, St.\ Louis, MO 63130.
    Email: {\scriptsize \texttt{\{r.he,hgonzale\}@wustl.edu}}.%
  }
}

\begin{document}

\maketitle
\thispagestyle{empty}
\pagestyle{empty}

\input{sections/abstract}
\input{sections/introduction}
\input{sections/problem_description}
\input{sections/tools}
\input{sections/examples}
\input{sections/conclusion}

\bibliographystyle{IEEEtran}
\bibliography{refs}

\end{document}

%% file: sections/abstract.tex
\begin{abstract}
  Efficiency, comfort, and convenience are three major aspects in the design of control systems for residential Heating, Ventilation, and Air Conditioning (HVAC) units.
  In this paper we propose an optimization-based algorithm for HVAC control that minimizes energy consumption while maintaining a desired temperature in a room.
  Our algorithm uses a Computer Fluid Dynamics model, mathematically formulated using Partial Differential Equations (PDEs), to describe the interactions between temperature, pressure, and air flow.
  Our model allows us to naturally formulate problems such as controlling the temperature of a small region of interest within a room, or to control the speed of the air flow at the vents, which are hard to describe using finite-dimensional Ordinary Partial Differential (ODE) models.
  Our results show that our algorithm produces significant energy savings without a decrease in comfort.
\end{abstract}

%% file: sections/introduction.tex
\section{Introduction}
\label{sec:intro}

Heating, Ventilation, and Air Conditioning (HVAC) systems are complex mechanical devices that control the climate of all kinds of buildings, large and small, residential and commercial.
In most situations, HVAC systems are used to maintain comfortable temperatures, while limiting both humidity and air speed away from undesirable levels, as described in standards such as ASHRAE~55~\cite{ASHRAE55}.
Yet, HVAC systems are typically controlled in a centralized or static fashion, disregarding variations in building configuration (windows or doors opened or closed dynamically), human activity, or even human perception of environmental conditions.
In this paper we present a control framework, based on a PDE-constrained optimal control problem, that takes into account localized conditions at the room or even person scale.
We also present simulated scenarios showing that our framework produces significant energy savings when compared to classical control strategies, and that it can react to changing occupation conditions.

There are many optimization-based studies of the control of HVAC units.
Goyal and Barooah~\cite{goyal2012method} studied in detail the use of RC network circuits to model the temperature within buildings.
Kelman and Borelli~\cite{kelman2011bilinear}, as well as Hazyuk et al.~\cite{hazyuk2012optimal,hazyuk2012optimal2}, used low-dimension ODE-based models to control a HVAC unit using Model Predictive Control (MPC).
Aswani et al.~\cite{aswani2012energy} used a learning-based MPC algorithm~\cite{aswani2013provably} to account for unmodeled dynamics and disturbance in ODE models when controlling HVAC units.
Domahidi et al.~\cite{domahidi2014learning} and Fux et al.~\cite{fux2014ekf} also used a learning-based method and MPC, the first using ADABOOST to estimate uncertainties and the second using an Extended Kalman Filter.
Ma et al.~\cite{ma2013stochastic} used Stochastic MPC to handle disturbances in the control of HVAC units, also using ODE models.
These results show that optimal control strategies increase the efficiency of HVAC systems, yet the use of ODE (i.e., concentrated parameter) models means that there is no detailed control of the air flow or of the temperature in arbitrary points in a room.

To address these issues, many authors have implemented optimization-methods using PDE (i.e., distributed parameter) models.
Even though the dynamic behavior of the air flow is mathematically complex due to turbulent dynamical responses, its response is laminar in larger areas~\cite{sinha2000numerical,bathe2004finite}, thus it can be analyzed using simpler non-turbulent Computer Fluid Dynamics (CFD) models.
Moreover, the existence and smoothness of the solutions of these non-turbulent models has been proved under suitable conditions~\cite{doering1995applied,foias2001navier}, thus they can be used together with gradient-based optimization algorithms.
Ito~\cite{ito1998} studied the theoretical optimal control of stationary Navier-Stokes equations coupled with the heat transfer equation, finding necessary conditions for the existence of an optimal argument.

Different variations of CFD models describing heating and ventilation situations in buildings can be found in the literature, such as the papers by Bathe et al.~\cite{bathe1995finite}, Sinha et al.~\cite{sinha2000numerical}, van Schijndel~\cite{van2011multiphysics}, Waring and Siegel~\cite{waring2008particle}, and the book by Awbi~\cite{awbi2003ventilation}.
Yet, these models are not suitable for integration with gradient-based optimization algorithms, which require the explicit formulation of all the approximating equations and their gradients.
It is for this reason that we formulated our own numerical discretization of the CFD model using the Finite Elements Method (FEM) and mixed boundary conditions.

In recent years, Burns et al.~\cite{borggaard2009,burns2012,burns2013} studied the optimal control of HVAC systems with linearized Navier-Stokes models around an arbitrary steady-state solution, and then used a Linear Quadratic Regulator (LQR) controller to find the optimal solution.
Their results show that PDE models, usually considered too complex for online numerical calculations, can be effectively used for building control.
Yet the authors' use of linearized approximations mean that the optimal solution is accurate only if it is close to the original steady-state linearization point, which poses a serious limitation in practical applications.
In this paper we use a nonlinear Navier-Stokes model, which means that our controller can produce large variations in the climate variables, at the cost of more expensive numerical computations.
As we show below, a careful choice of parameters and numerical algorithms allows for the nonlinear numerical computations to be performed within the time constant of this dynamical system.

This paper presents two contributions.
First, we use a nonlinear non-turbulent Navier-Stokes model together with a convection-diffusion heat equation to model the climate in a building.
We discretize this model using a Finite Element Method (FEM), and we use it to find the optimal control for the HVAC system in a building.
Second, we show via simulations that it is possible to design control objectives that are functions of the residents' locations, which greatly increase the efficiency of the HVAC system.


The paper is organized as follows: Section~\ref{sec:probdesc} describes the CFD model and the formulation of the optimal control problem; Section~\ref{sec:tools} discusses the discretization of the PDEs forming the CFD model, and the discretization of the optimal control problem; and Section~\ref{sec:examples} gives the results of our simulated experiments after we incorporate small changes in the actuators of the HVAC unit, such as independent control of each vent in a room or control of the angle of the air flow in a vent.
Our results validate our hypothesis that sizable energy savings can be obtained by introducing small improvements in the actuation of HVAC units, thanks to the accurate CFD models that describe the dynamical and distributed behavior of the climate variables in a building.

%% file: sections/problem_description.tex
\section{Problem Description}
\label{sec:probdesc}

A commonly missing key feature in many physical climate models used to control HVAC systems is the ability to capture the real-time spatial variability of the temperature and air flow, depending on the floor plan and configuration of the building (e.g., open or closed door and windows).
For this reason, we use a CFD model, which explicitly considers temporal and spacial variations, to describe the interactions between the temperature, air flow, and pressure.
We then formulate an optimal control problem where our CFD model appears as a constraint, and whose objective function aims to minimize the energy consumption of the HVAC system while maintaining the temperature constant at a desired reference.
The inclusion of spatial variables in our description of the climate variables will not only improve the accuracy of our estimations, but will also allow us to naturally formulate richer problems, such as focusing only on an specific region in a room, as shown in Section~\ref{sec:examples}.

In the remainder of this section we introduce in detail the CFD model and the optimal control problem.

\subsection{CFD Model}
\label{sec:pde_model}

The foundation of our model is the Navier-Stokes equation, which couples temperature with free flow convection~(as explained Section~8 in~\cite{awbi2003ventilation}, among other references).
As shown in the literature, atmospheric air can be modeled as an incompressible Newtonian fluid when the temperature is between $-20\,\unit{^\circ C}$ and $100\,\unit{^\circ C}$~\cite{awbi1989,Dobrzynski2004}.
Hence, we can use the Navier-Stokes equation for incompressible laminar flows, together with the convection-diffusion temperature model for fluids.

Throughout the paper we make two major simplifications to the CFD model.
First, we assume that the air flow behaves as a laminar fluid which reaches a steady-state behavior much faster than the temperature in the building.
As mentioned in Section~\ref{sec:intro}, both laminar and turbulent flows are present in general in a residential building, for example, in the area around HVAC vents~\cite{awbi2003ventilation}.
However, Sun et al.~\cite{sun2003} found only minor differences between laminar and turbulent models in a geometry similar to ours, while turbulent models are significantly more complex than laminar models~\cite{akhtar2010high,clarke2001energy}.
Hence, we consider a stationary Navier-Stokes equation to describe the fluid behavior, and a time-dependent equation to describe the temperature behavior.
Second, we consider only two-dimensional air flows moving parallel to the ground.
This assumption intuitively makes sense since the air flow in the top half of a room can be accurately estimated using a two-dimensional model, mostly due to the lack of obstacles (such as furniture).
These assumptions reduce the accuracy of our model to some extent, e.g., van der Poel et al.~\cite{van2013comparison} compared 2D and 3D Rayleigh-B\'{e}nard convection simulations for a cylindrical geometry and showed that differences arise for Prandtl constants $\prno < 1$, while our model's Prandtl constant is $\prno = 1.2$.
Yet, both assumptions allow us to significantly simplify the computational complexity of our CFD-based control design (measured by the number of variables and number of equality constraints of the model), which in turn allows us to compute results on the order of tens of minutes, as shown in Section~\ref{sec:examples}.

Let $\Omega \subset \R^2$ be the area of interest.
We will denote the \emph{boundary of $\Omega$} by $\bdry{\Omega}$.
Let $u\colon \Omega \to \R^2$ be the \emph{stationary air flow velocity}, and $p\colon \Omega \to \R$ be the \emph{stationary air pressure} in $\Omega$.
Also, given $T > 0$, let $T_e\colon \Omega \times [0,T] \to \R$ be the \emph{temperature} in $\Omega$.
Then, using the formulation found in~\cite{landau_1971}, the convection-diffusion of temperature in $\Omega$ can be described by the following PDE:
\begin{multline}
  \label{eq:temp}
  \pderiv{T_e}{t}(x,t) - \diver[x]\pb{\kappa(x)\, \grad[x]{T_e}(x,t)} + u(x) \cdot \grad[x]{T_e}(x,t) =\\
  = g_{T_e}(x,t),
\end{multline}
where $g_{T_e}\colon \Omega \times [0,T] \to \R$ represents the heat sources in the room, $\kappa\colon \Omega \to \R$ is the \emph{thermal diffusivity}, $\diver[x]{} = \pderiv{}{x_1} + \pderiv{}{x_2}$ is the \emph{divergence operator}, and $\grad[x]{} = \pb{\pderiv{}{x_1},\pderiv{}{x_2}}^T$ is the \emph{gradient operator}.

Similarly, the stationary air flow in $\Omega$ is governed by the following set of incompressible Navier-Stokes stationary PDEs:
\begin{multline}
  \label{eq:ns1}
  - \frac{1}{\reno}\, \lapl[x]{u}(x) + \pb{u(x) \cdot \grad[x]{}}\, u(x) + \frac{1}{\rho}\, \grad[x]{p(x)} + \alpha(x)\, u(x) =\\
  = g_u(x);\ \text{and},
\end{multline}
\begin{equation}
  \label{eq:ns2}
  \diver[x]{u}(x) = 0,
\end{equation}
where $g_u\colon \Omega \to \R^2$ represents all the external forces applied to the air (such as fans), $\reno$ is the \emph{Reynolds number} (which is inversely proportional to the \emph{kinematic viscosity}), $\rho$ is the \emph{density of the air}, $\alpha\colon \Omega \to \R$ is the \emph{friction constant}, $u(x) \cdot \grad[x]{} = u_1(x)\, \pderiv{}{x_1} + u_2(x)\, \pderiv{}{x_2}$ is the \emph{advection operator}, and $\lapl[x]{} = \pderiv[2]{}{x_1} + \pderiv[2]{}{x_2}$ is the \emph{Laplacian operator}.
Since we do not model the vertical dimension of $\Omega$, we omit the Boussinesq-approximation buoyancy term proportional to $T_e$, which is typically included on the right-hand side of~\eqref{eq:ns1}.
We modify $\kappa$ and $\alpha$ to model obstacles to heat and air flow in $\Omega$, such as walls, doors, and windows, as described in~\cite{gersborg2005,pingen2007topology}.
In particular, when the point $x$ corresponds to a material that blocks air, we choose $\alpha(x) \gg u(x)$, which results in $u(x) \approx 0$, and when the point $x$ corresponds to air then we choose $\alpha(x) = 0$.

We divide the boundary of $\Omega$ to two outlets of the HVAC system, denoted by $\Gamma_o$, one air return inlet, denoted by $\Gamma_i$, and the exterior walls, denoted by $\Gamma_w$.
Thus $\Gamma_i \cup \Gamma_o \cup \Gamma_w = \bdry\Omega$.

We use a mix of boundary conditions to model the effect of the HVAC system in the room, as explained below.
Let $\hat{n}(x)$ be the inward-pointing unit vector perpendicular to the boundary at $x \in \bdry{\Omega}$.
Hence, the air flow has the following boundary conditions:
\begin{itemize}
\item The HVAC unit's fan sets the air flow at $\Gamma_o$, hence $u(x) = u_o\, \hat{n}(x)$ for each $x \in \Gamma_o$, where $u_o > 0$ is the HVAC fan speed.
\item The airflow at the inlet is not constrained, hence $u(x)$ is free for each $x \in \Gamma_i$.
\item The airflow satisfies a no-transverse condition at the exterior walls, hence $u(x) \cdot \hat{n}(x) = 0$ for each $x \in \Gamma_w$.
\end{itemize}
The boundary condition for the temperature is $T_e(x) = T_A$ for each $x \in \bdry\Omega$, where $T_A$ is the \emph{atmospheric temperature}.
We apply a boundary condition for the pressure equation only at the inlet, setting $p(x) = p_A$ for each $x \in \Gamma_i$, where $p_A$ is the \emph{atmospheric pressure}.

\subsection{Optimal Control Problem}
\label{sec:ocp}

To control the temperature in the zone $\Omega_z \subset \Omega$, we aim to minimize the following cost function:
\begin{equation}
  \label{eq:cost_function}
  \int_0^{t_f} \!\p{\int_{\Omega_z} \pb{T_e(x,t) - T_e^*}^2 \diff{x} + \lambda_1\, v^2(t)} \diff{t} + \lambda_2\, u^2_o,
\end{equation}
where $\lambda_{1,2} > 0$, $T_e^*$ is the reference temperature set by the user, and $t_f$ is the time horizon.
The heater power $v(t)$ and the fan speed $u_o$ are our controlled variables, the former appearing as $g_{T_e}(x,t) = v(t)$ for each $x \in \Theta_h \subset \Omega$, and the latter appearing as a boundary condition.
We formulate a PDE-constrained optimal control problem using the cost in~\eqref{eq:cost_function}, together with the CFD model in~\eqref{eq:temp}-\eqref{eq:ns2} and its boundary conditions as constraints.
We also add inequality box constraints for all the controlled variables, so they remain within safety limits.

As explained in Section~\ref{sec:examples}, our experiments introduce variations to the cost function in~\eqref{eq:cost_function} depending on the number of available actuators.
Regardless, the goal of regulating the temperature will remain the same throughout all our experiments.

%% file: sections/tools.tex
\section{Numerical Implementation}
\label{sec:tools}

Our numerical implementation of the PDE-constrained optimal control problem described in Section~\ref{sec:ocp} is obtained by first using FEM to transform the CFD model in~\eqref{eq:temp}-\eqref{eq:ns2} to a set of ODEs as described in Chapters~3 and~4 of~\cite{cuvelier_book}, and then using the consistent approximation technique described in Chapter~4 of~\cite{polak1997} which transforms optimal control problems (with ODE constraints) into nonlinear programming problems.
After those two transformations, we use commercially available numerical solvers to find approximations of the desired optimal control, as described in Section~\ref{sec:examples}.

\subsection{FEM Discretization}

Among the many discretization techniques for PDEs, FEM stands out for being compatible with complex geometries of the domain $\Omega$.
Intuitively speaking, FEM approximates PDEs by dividing the domain into polygons, and then finding a set of ODEs for each vertex, and possibly each facet, of each polygon.
The resulting set of ODEs has the property that each ODE is dependent only on its neighbors.

Before we can formally describe the FEM discretization, we need to introduce extra notation.
Let $H^1(\Omega,\R^n)$ be the set of functions from $\Omega$ to $\R^n$ belonging to $L^2(\Omega,\R^n)$, whose weak derivative is also in $L^2(\Omega,\R^n)$~\cite{Ziemer1989}.
Note that $H^1(\Omega,\R^n)$, endowed with the dot product $\dprod{f}{g} = \int_\Omega f(x) \cdot g(x)\, \diff{x}$, is a Hilbert space.
Similarly, we denote $\dprod{f}{g}_S = \int_S f(x) \cdot g(x)\, \diff{x}$.

Let $\set{W_k}_{k=1}^{N_{pl}}$ be a polygonal partition of $\Omega$, i.e., $\bigcup_{k=1}^{N_{pl}} W_k = \Omega$, $\sint(W_k) \cap \sint(W_j) = \emptyset$ for each $k \neq j$, and each $W_k$ is a polygon.
If $\set{x_k}_{k=1}^{N_v}$ is the set of vertices in the polygonal partition and $\set{y_j}_{j=1}^{N_w}$ is the set of nodal points, then we define the \emph{test functions} $\set{\xi_k}_{k=1}^{N_v}, \set{\psi_k}_{k=1}^{N_v} \subset H^1(\Omega,\R)$, and $\set{\phi_k}_{k=1}^{2N_w} \subset H^1(\Omega,\R^2)$, where $N_v,N_w \in \N$ and $N_v \leq N_w$, with the following properties for each $k \in \set{1,\dotsc,N_v}$ and each $j \in \set{1,\dotsc,N_w}$:
\begin{itemize}
\item $\xi_k$, $\phi_k$, and $\psi_k$ are continuous;
\item $\xi_k$, $\phi_k$, and $\psi_k$ are nonzero only in the polygons containing $x_k$; and,
\item $\xi_k(x_k) = \psi_k(x_k) = 1$, $\phi_{2j-1}(y_j) = \smats{1\\0}$, and $\phi_{2j}(y_j) = \smats{0\\1}$.
\end{itemize}
Then, from~\eqref{eq:temp}-\eqref{eq:ns2}, and using Green's Formulas (see Appendix~C.2 in~\cite{Evans2010}), we get the following Galerkin identities (as described in Chapter~3.6 of~\cite{cuvelier_book}):
\begin{multline}
  \label{eq:temp_fem}
  \dprod{\pderiv{T_e}{t}(\cdot,t)}{\xi_k}
  - \dprodb{\kappa(x)\, \grad[x]{T_e}(\cdot,t)}{\grad[x]{\xi_k}} +\\
  + \dprodb{u \cdot \grad[x]{T_e}(\cdot,t)}{\xi_k}
  = \dprodb{g_{T_e}(\cdot,t)}{\xi_k};
\end{multline}
\begin{multline}
  \label{eq:ns1_fem}
    -\frac{1}{\reno}\, \dprod{\grad[x]{u}}{\grad[x]{\phi_j}}
    + \dprod{\pb{u \cdot \grad[x]{}}\, u}{\phi_j}
    + \dprod{\grad[x]{p}}{\phi_j} +\\
    + \dprod{\alpha\, u}{\phi_j}
    = \dprod{g_u}{\phi_j};\ \text{and},
\end{multline}
\begin{equation}
  \label{eq:ns2_fem}
  \dprod{\diver[x]{u}}{\psi_k} = 0,
\end{equation}
for each $k \in \set{1,\dotsc,N_v}$, $j \in \set{1,\dotsc,2\, N_w}$, and almost every $t \in \sqparen{0,t_f}$.

Now, given $N_{T_e}, N_{u}, N_{p} \in \N$, consider the linearly independent sets of \emph{basis functions} $\setb{\bfn{\xi}_j}_{j=1}^{N_{T_e}}, \setb{\bfn{\psi}_j}_{j=1}^{N_u} \subset H^1(\Omega,\R)$, and $\setb{\bfn{\phi}_j}_{j=1}^{N_p} \subset H^1(\Omega,\R^2)$.
Using these basis functions we can project the variables of our CFD model into finite-dimensional subspaces, i.e.:
\begin{multline}
  \label{eq:fdim_reps}
  T_e(x,t) = \sum_{j=1}^{N_{T_e}} \eta_{T_e,j}(t)\, \bfn{\xi}_j(x), \
  u(x) = \sum_{j=1}^{N_u} \eta_{u,j}\, \bfn{\phi}_j(x),\\
  p(x) = \sum_{j=1}^{N_p} \eta_{p,j}\, \bfn{\psi}_j(x).
\end{multline}
Applying the representations in~\eqref{eq:fdim_reps} to the Galerkin identities in~\eqref{eq:temp_fem} results in a set of $N_v$ ODEs with state variables $\set{\eta_{T_e,j}}_{j=1}^{N_{T_e}}$.
Similarly, applying the representations to~\eqref{eq:ns1_fem}-\eqref{eq:ns2_fem} results in a set of $2\, N_w$ nonlinear algebraic equations with parameters $\set{\eta_{u,j}}_{j=1}^{N_u}$ and $N_v$ linear ones with parameters $\set{\eta_{p,j}}_{j=1}^{N_p}$.
All these differential and algebraic equations are, in practice, parametrized by constants corresponding to the inner products between basis and test functions, as well as their gradients.
We omit the technical details of the final set of equations due to space constraints, and we refer the interested reader to Chapter~3 in~\cite{cuvelier_book} for more information.

\subsection{Optimal Control Discretization}

After the FEM discretization, we effectively have an Differential Algebraic Equation (DAE) optimal control problem where~\eqref{eq:temp_fem} contributes $N_v$ ODEs, \eqref{eq:ns1_fem} contributes $2\, N_w$ nonlinear equality constraints, and~\eqref{eq:ns2_fem} contributes $N_v$ linear constraints.
Several extra equality constraints are added due to the boundary conditions of the air flow and the pressure, as described in Section~\ref{sec:pde_model}.
The actual number of constraints due to boundary conditions depends on the number of vertices in the polygonal partition $\set{W_k}_{k=1}^{N_{pl}}$ over the boundary.

The consistent approximation of this type of optimal control problem is studied in Chapter~4 of~\cite{polak1997}.
We follow the procedure described there, i.e., we first normalize the problem using the technique described in Chapter~4.1.2 of the same book, and then we use the Forward-Euler discretization method to transform the ODEs into a sequence of equality constraints.
Again, we omit the technical details of the final equality-constrained nonlinear programming problem due to space constraints.

%% file: sections/examples.tex
\section{Experimental Results}
\label{sec:examples}

We simulated a two-room apartment with a square area of interest~$\Omega_z$ (e.g., the area were a resident is located).
Our goal is to show that using zoned control over the area of interest produces a significant improvement over controlling the temperature over the whole room.
The HVAC system consists of two heaters and two forced-air outlets with variable-speed fans.
Moreover, in order to show the efficiency and stability of our zoned control algorithm, we simulated 18~different scenarios with different areas of interest, distributed uniformly over the apartment.

\begin{figure}[tp]
  \centering
  \includegraphics[width=.9\linewidth,trim=120 160 100 160,clip]{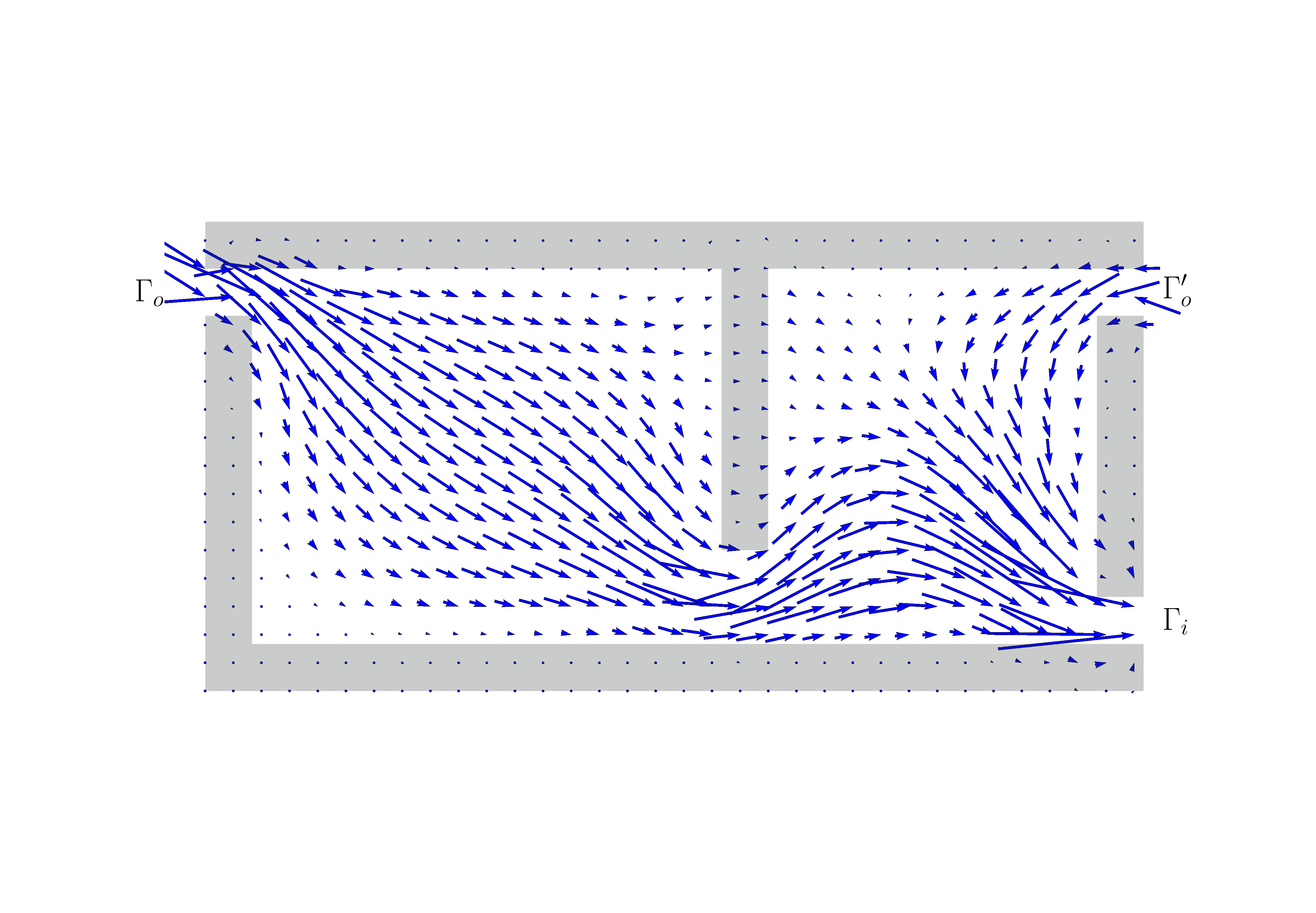}%
  \caption{%
    Air flow of the experiment in Section~\ref{sec:exp1} with target area $\Omega$, i.e., the whole apartment.
    Average air speed in the apartment is~$0.11 \unit{m/s}$.
    Walls are shown in shaded black.
    Two outlets and one inlet of the HVAC system are marked as $\Gamma_o$,
    $\Gamma_o'$ and $\Gamma_i$ respectively.
  }
  \label{fig:exp1_airflow}
\end{figure}

A diagram of the apartment is shown in Figures~\ref{fig:exp1_temp} and~\ref{fig:exp1_airflow}.
The apartment's dimensions are $5 \times 10 \unit{m^2}$, the width of all outlets and inlet is $0.5 \unit{m}$, and each of the areas has dimensions $2 \times 2 \unit{m^2}$.
The two heaters are denoted by $\Theta_h$ (left) and $\Theta_h'$ (right), with dimensions $1 \times 1 \unit{m^2}$.

The fluid mechanics are governed by the constants $\reno = 0.05$ and $\prno = 1.2$, with $k(x) = 10^{-2}$ and $\alpha(x) = 0$ when $x \in \Omega$ is located in free air, while $k(x) = 10^{-4}$ and $\alpha(x) = 100$ when $x \in \Omega$ is located on or in a wall.
The atmospheric pressure is $p_A = 101.3 \unit{kPa}$, and the atmospheric temperature is $T_A = 23.83 \unit{^\circ C}$.
We set the desired temperature to $T_e^* = 24.83 \unit{^\circ C}$, the time horizon to $t_f = 300 \unit{s}$, and we assume the function $g_u$ is identically zero (i.e., no fans are inside the room).
As explained in Section~\ref{sec:ocp}, $g_{T_e}(x,t) = v(t)$ for each $x \in \Theta_h$, $g_{T_e}(x,t) = v'(t)$ for each $x \in \Theta_h'$, and $g(x,t) = 0$ otherwise.
The parameters in the cost function~\eqref{eq:cost_function} are $\lambda_1 = 0.002$ and $\lambda_2 = 0.001$.
The optimal control problem finds the fan speeds $u_o$ and $u_o'$ for $\Gamma_o$ and $\Gamma_o'$, respectively, and the heater powers $v(t)$ and $v'(t)$ for $\Theta_h$ and $\Theta_h'$, respectively.
We set the fan speed box constraints to $\sqparenb{0.1,1} \unit{\frac{m}{s}}$, and the heater power box constraints to $\sqparenb{0,5} \unit{kW}$.

We discretized the area into $N_{pl} = 452$ elements, and the number of total nodes is $N_v = 227$.
We used first-order Lagrangian elements to define the test and basis functions $\xi_k$, $\bfn{\xi}_k$, $\psi_k$, and $\bfn{\psi}_k$, thus $N_{T_e} = N_p = 227$.
We used second-order Lagrange elements to define the test and basis functions $\phi_k$ and $\bfn{\phi}_k$, thus $2\, N_w = N_u = 1696$.
More details regarding our choice of test and basis functions can be found in Chapter~3.3.1 of~\cite{Logg2012}.
The ODE discretization time step was chosen as $\Delta t = 10 \unit{s}$.

We calculated the total energy usage as the sum of the heater energy usage, i.e., $\int_0^{t_f} v(t) + v'(t)\, \diff{t}$, and the fan energy usage as $\int_0^{t_f} \int_{\Gamma_o \cup \Gamma_o'} \norm{u(x)}\, p(x)\, \diff{x}\, \diff{t}$.
Our results were obtained using a 16-core \emph{Xeon E5-2680} computer running at $2.7 \unit{GHz}$, with $128 \unit{GB}$ of RAM.
We wrote our code using Python, the FEM discretization was computed using tools from the \emph{FEniCS Project}~\cite{Logg2012}, and the nonlinear programming problem was numerically solved using the \emph{SNOPT} library~\cite{gill_2002} interfaced using the \emph{OptWrapper} library~\cite{optwrapper}.
The computation time ranged between $15 \unit{min}$ and $45 \unit{min}$ for each experiment.

\newcommand{\scb}{0.32}
\newcommand{\scc}{0.47}
\begin{figure*}[tp]
  \centering
  \subfloat[%
  Temperature distribution at time $t_f$ with target area $\Omega$, i.e., the whole apartment.
  ]{%
    \label{fig:exp1_temp_whole}%
    \includegraphics[width=\scb\linewidth,trim=0 0 30 0,clip]{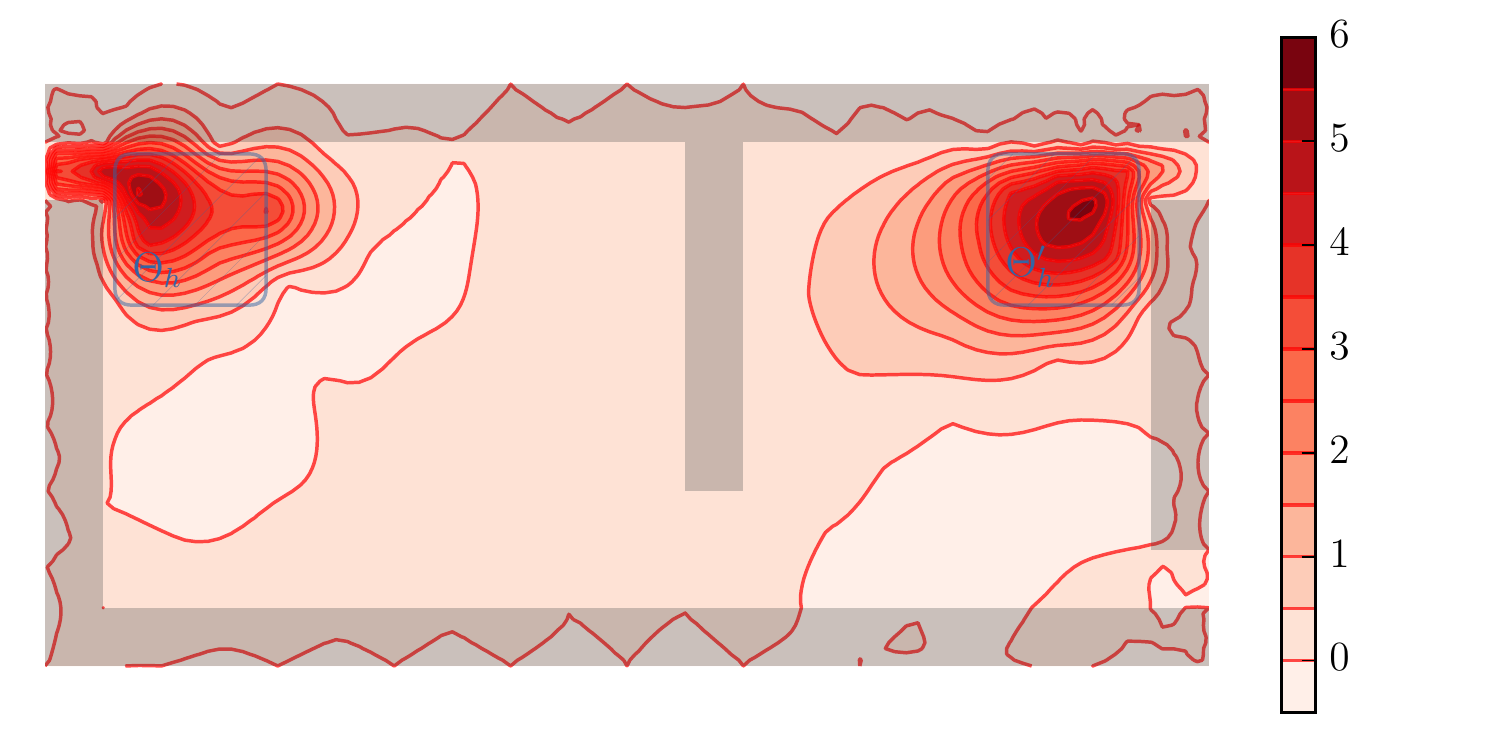}%
  }%
  \hfill%
  \subfloat[%
  Temperature distribution at time $t_f$ with target area $\Omega_z$.%
  ]{%
    \label{fig:exp1_temp_02}%
    \includegraphics[width=\scb\linewidth,trim=0 0 30 0,clip]{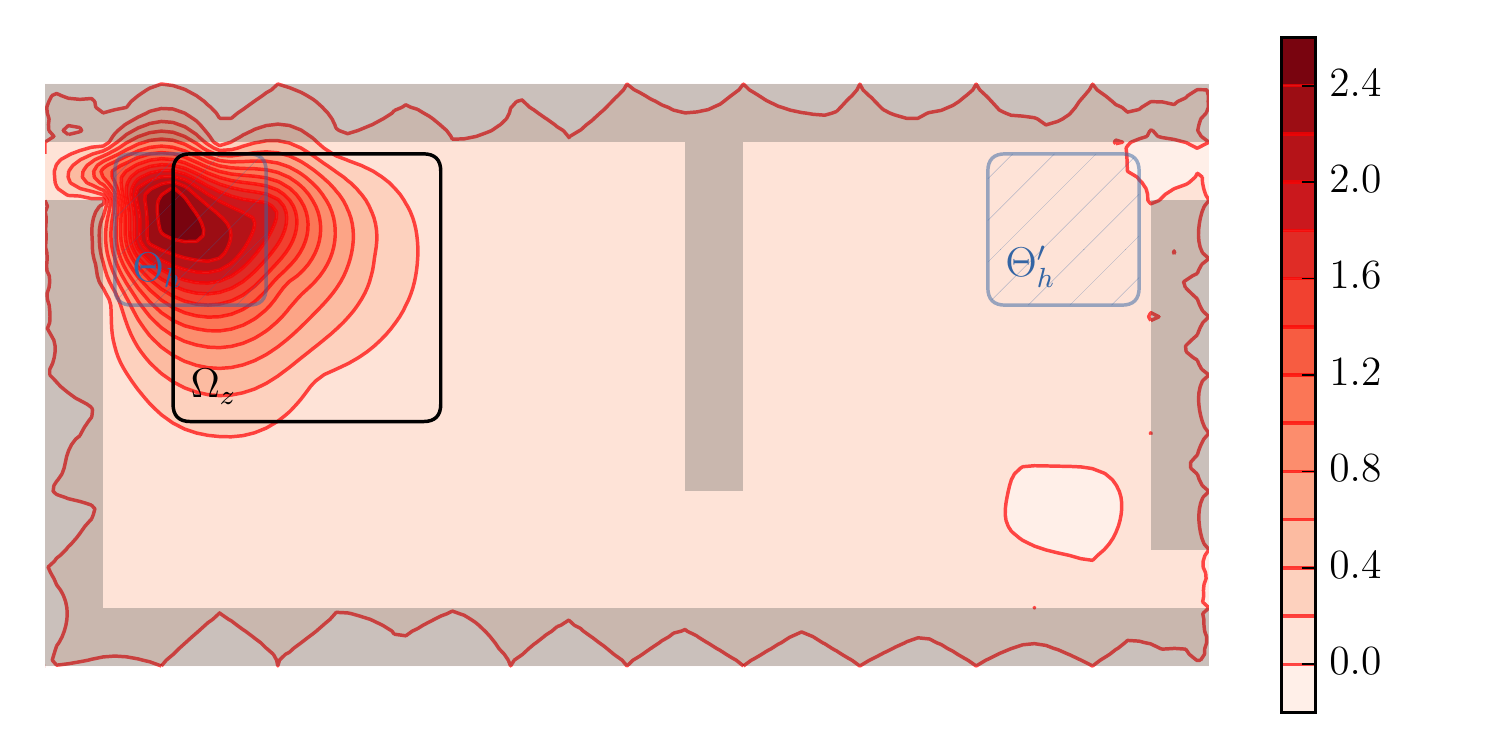}%
  }%
  \hfill%
  \subfloat[%
  Temperature distribution at time $t_f$ with target area $\Omega_z$.%
  ]{%
    \label{fig:exp1_temp_07}%
    \includegraphics[width=\scb\linewidth,trim=0 0 30 0,clip]{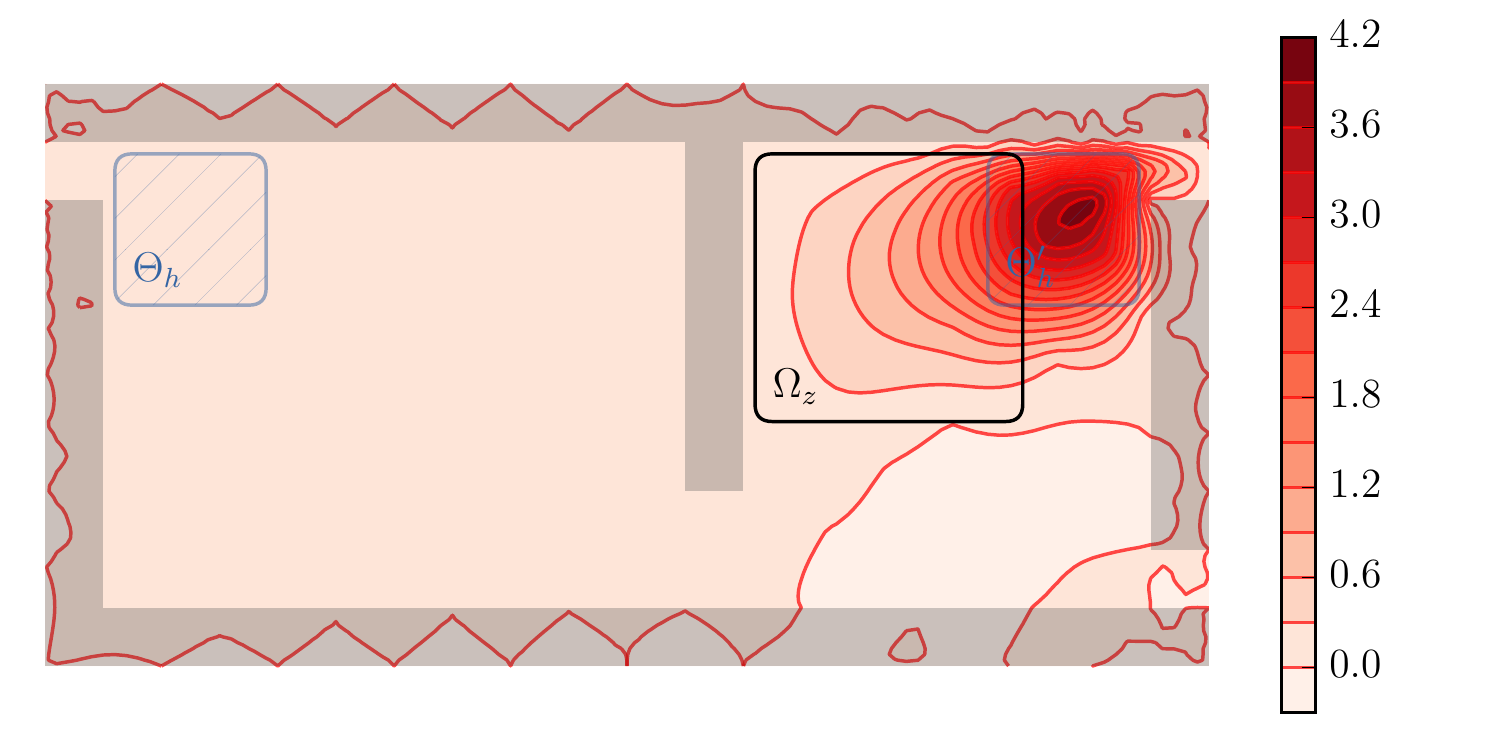}%
  }%
  \caption{%
    Results of the experiments in Section~\ref{sec:exp1}.
    Walls are shown in shaded black, and heaters are shown in shaded blue.
    Values are in $\unit{^\circ C}$ with respect to $T_A$.
  }
  \label{fig:exp1_temp}
\end{figure*}

\subsection{Nonlinear Navier-Stokes Model}
\label{sec:exp1}

Using the nonlinear Navier-Stokes model described in \eqref{eq:temp}-\eqref{eq:ns2} we simulated two major scenarios, the first where the objective function in~\eqref{eq:cost_function} uses $\Omega_z = \Omega$, i.e., the target area is the whole apartment, and the second where we use a zoned approach with a smaller $\Omega_z$ which moves around the apartment to 18~different locations, as explained above.
Figure~\ref{fig:exp1_temp_whole} shows the temperature distribution for the first scenario, while Figures~\ref{fig:exp1_temp_02} and~\ref{fig:exp1_temp_07} show the temperature distribution for two of the 18~zoned simulations.
Also, Figure~\ref{fig:exp1_airflow} shows the stationary airflow for the first scenario.

\begin{figure}[tp]
  \centering
  \subfloat[%
  Average absolute temperature error, in $\unit{^\circ C}$, within $\Omega_z$ at time $t_f$ with respect to $T_e^*$.
  ]{%
    \label{fig:exp_stats_fin}%
    \includegraphics[width=\scc\linewidth,trim=0 10 0 10,clip]{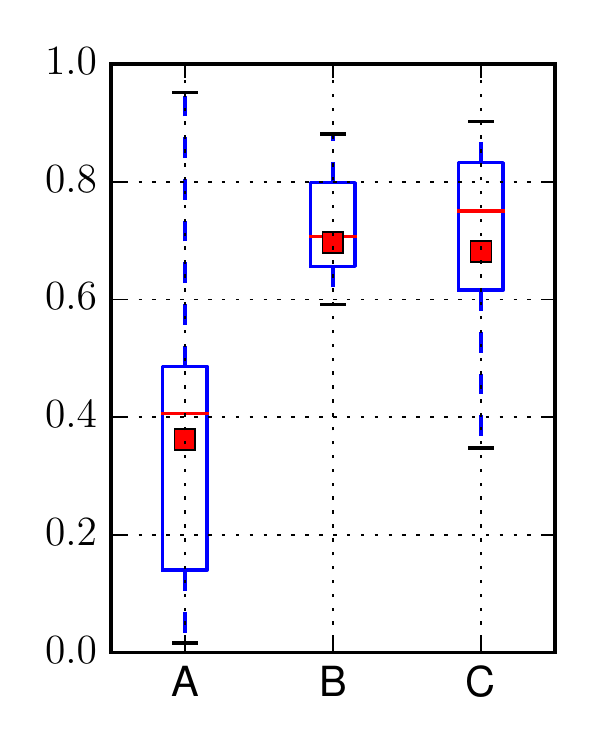}%
  }%
  \hfill%
  \subfloat[%
  Ratio of energy usage over average temperature change within $\Omega_z$,
  in $\unit{Wh/^\circ C}$.
  ]{%
    \label{fig:exp_stats_eng}%
    \includegraphics[width=\scc\linewidth,trim=0 10 0 10,clip]{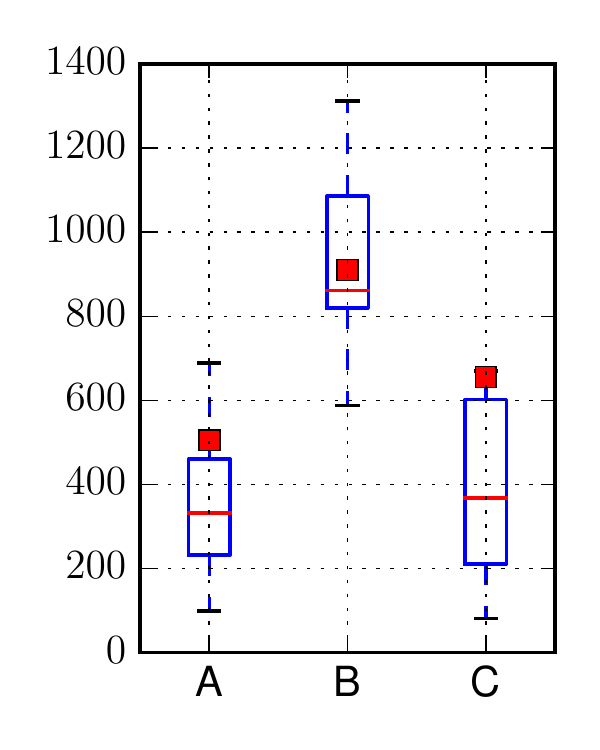}%
  }%
  \caption{%
    Results of the experiments in Section~\ref{sec:exp1} and~\ref{sec:exp2}.
    Columns: (A)~nonlinear model with zoned control, (B)~nonlinear model without zoned control, (C)~linearized model with zoned control.
    Each column shows the median (red line), mean (red box), and first-to-third quartiles (blue box).
  }
  \label{fig:exp_stats}
\end{figure}

In the first scenario, where $\Omega_z = \Omega$, the optimal average absolute temperature error in the apartment was $0.502 \unit{^\circ C}$ at time $t_f$, and the ratio of energy usage over average temperature change within $\Omega$ was $1009.2 \unit{Wh/^\circ C}$.
We calculated the same statistics for the 18~different zones $\Omega_z$, which are summarized in Column~(A) of Figures~\ref{fig:exp_stats_fin} and~\ref{fig:exp_stats_eng}.
To make both scenarios comparable, we recalculated these statistics for the first scenario, this time considering the average temperature changes in $\Omega_z$ instead of $\Omega$, which are summarized in Column~(B) of Figures~\ref{fig:exp_stats_fin} and~\ref{fig:exp_stats_eng}.
Those figures clearly show that using zoned control is significantly more accurate and more efficient than heating the whole apartment.
It is worth noting that the zoned approach requires roughly half the energy to change the average temperature by $1 \unit{^\circ C}$ in $\Omega_z$ when compared to the first scenario.

The results in Figure~\ref{fig:exp1_temp} indicate that when the resident is near one heater, say $\Theta_h$, our algorithm automatically shut down the other heater, say $\Theta_{h'}$, as intuitively expected.
Therefore, if it is possible to localize a resident within an apartment, e.g., via Bluetooth beacons or using a sensor network, then we can increase the efficiency of the HVAC unit significantly without major modifications to the mechanical ventilation system.

\subsection{Linearized Navier-Stokes Model}
\label{sec:exp2}

We also computed the optimal control using the linearized Navier-Stokes model described in~\cite{burns2012,burns2013}.
The linearized model has clear advantages over our nonlinear model, including a larger set of theoretical results supporting it, and a faster computation time.
On the other hand, linearized models perform well only when the values produced by the model are close to the stationary linearization point.

\begin{figure}[tp]
  \centering
  \includegraphics[width=.8\linewidth,trim=0 0 30 0,clip]{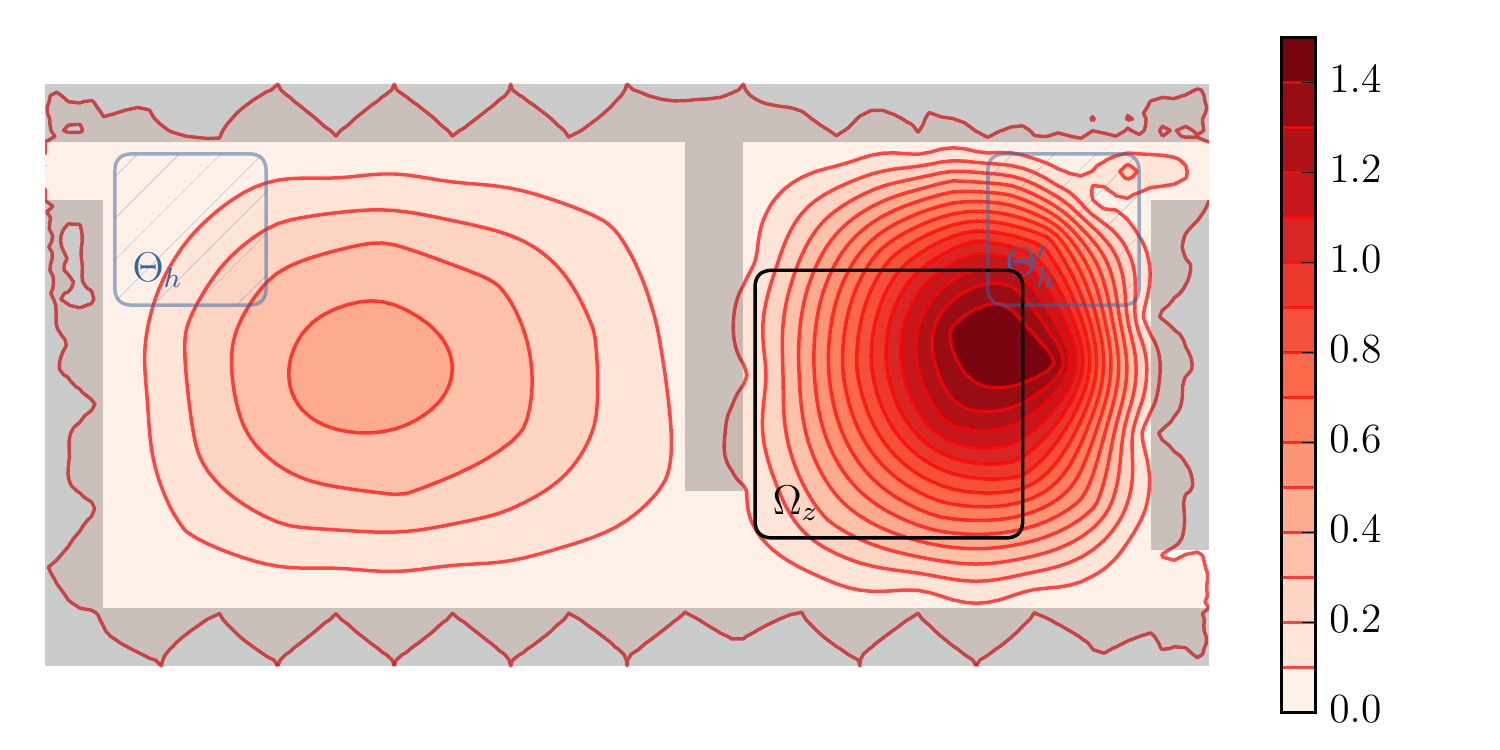}%
  \caption{%
    Absolute error in temperature distribution, with respect to \emph{FEniCS} simulation, for the linearized model in Section~\ref{sec:exp2}, at time $t_f$ with target area~$\Omega_z$.
    Values are in $\unit{^\circ C}$.%
  }%
  \label{fig:exp2_err}
\end{figure}

We ran the same experiments as in the second scenario of Section~\ref{sec:exp1}, i.e., controlling the temperature in 18~different zones.
The statistics for average absolute temperature error in $\Omega_z$, and ratio of energy usage over average temperature in $\Omega_z$, are shown in Column~(C) of Figures~\ref{fig:exp_stats_fin} and~\ref{fig:exp_stats_eng}.
Even though the energy efficiency is comparable when we use linearized or nonlinear models, the accuracy is significantly different, with the nonlinear model consistently performing better than the linearized model.
We believe the difference is due to the lack of accuracy of the linearized model.
As exemplified in Figure~\ref{fig:exp2_err}, there is a large error in temperature distribution between the linearized model and an accurate benchmark simulation using the \emph{FEniCS} solver, which is consistent with the highly nonlinear behavior of the Navier-Stokes equation.

%% file: sections/conclusion.tex
\section{Conclusion}
\label{sec:conclusion}

Our results open the door to a large number of exiting opportunities to improve the energy efficiency of buildings.
By making small improvements to existing HVAC units it is possible to dramatically increase the efficiency of HVAC units without a decrease in human comfort.
It is worth noting that our results do not require, in principle, the use of expensive variable-speed fans or variable-power heaters, since those control signals can be implemented using switched strategies~\cite{Vasudevan2013a,Vasudevan2013b}.
More importantly, our simplifications have allowed us to obtain results in tens of minutes while still capturing the distributed behavior of the climate variables, which is much closer to real-time applications than previous results~\cite{Dobrzynski2004}, achieving a good trade-off when compared to less accurate linearized Navier-Stokes models.